\documentstyle[12pt]{article}
\newcommand{\R}{I\!\!R}
\newcommand{\N}{I\!\!N}

\def\bgneqy{\begin{eqnarray}}
\def\endeqy{\end{eqnarray}}
\def\bgneqy*{\begin{eqnarray*}}
\def\endeqy*{\end{eqnarray*}}
\newcommand{\C}{I\!\!\!\!C}

\newcommand{\p}{\partial}

\newtheorem{thm}{Theora}[section]
\newtheorem{Th}[thm]{\normalsize{\bf{Theorem}}}

\newtheorem{Lem}[thm]{\normalsize{\bf{Lemma}}}
\newtheorem{Prop}[thm]{\normalsize{\bf{Proposition}}}
\newtheorem{Def}[thm]{\normalsize{\bf{Definition}}}

\newtheorem{Rem}[thm]{Remark}

\newcommand\dint{\displaystyle\int}

\newcommand\dfrac{\displaystyle\frac}

\begin{document}
%\articletype{}
\title{\textbf{ Uncertainty principle for the Poly-axially $L_{\alpha}^{2}$-multiplier operators}
\author{{Belgacem Selmi \textsuperscript{ab} \thanks{ Email:
Belgacem.Selmi@fsb.rnu.tn},\quad Rahma Chbeb\textsuperscript{ b}
\thanks{ Email: rahmachbeb@yahoo.fr}}
\date{ \small{
\textsuperscript{a}   University of Carthage, Faculty of Sciences of
Bizerte, 7021 Zarzouna, Departement of Mathematics.\\
\textsuperscript{b} University of Tunis El Manar, Faculty of Sciences of Tunis, Research Laboratories
of Mathematics Analysis and Applications LR11ES11, Tunisia.}}}}
%\date{}
\maketitle

\begin{abstract}
In this article, the Heisenberg-Pauli-Weyl uncertainty
principle and Donoho-Stark's uncertainty principle are obtained for the Poly-axially
$L^{2}$ multiplier operators.
\end{abstract}
\indent
%{\it{\bf {2010 AMS Mathematics Subject Classification:}}}\quad  46E35.\\[1mm]\indent
%%%%%%%%%%%%%%%%%%%%%%%%%%%%%%%%%%%%%%%%%%%%%%%%%%%%%%%%%%%%%%%%%%%%%%%%%%%%%%%%%%%%%%%%%%%%%%%%%%%%%%%%%%%%%%%%%
\section{Introduction}
\quad\quad Uncertainty principles are an essential restriction in
Fourier analysis. They state that a function $f$ and its Fourier
transform $\widehat{f}$ cannot be at the same time simultaneously
and sharply localized. That is, it is impossible for a non-zero
function
and its Fourier transform to be simultaneously small. Many mathematical formulations of this general fact can be
found in \cite{folland,havin,price}.\\
The classical uncertainty principle was established by Heisenberg \cite{hei}
bringing a fundamental problem in quantum mechanics. The
mathematical formulation was established by Kennard \cite{ken} and
Weyl (who attributes results to Pauli) \cite{weyl} in the late
1920's. This leads to this quantitative formulation in the form of a
lower bound of the product of the dispersions of a nonzero function
$f$ and its Fourier transform $\widehat{f}$:
%\begin{equation}\label{unc}
$$\Big\|\|x\|\,f\Big\|_{L^{2}(\R^{n})}\Big\|\|y\|\,\widehat{f}\Big\|_{L^{2}(\R^{n})}\geq\,
\frac{n}{2}\,\|f\|_{L^{2}(\R^{n})}^{2},$$%\, f\in L^{2}(\R^{n}).
%\end{equation}
Where $$\widehat{f}(\lambda)=(2\pi)^{\frac{-n}{2}}\int_{\R^{n}}f(x)e^{-i\langle x,\lambda\rangle}dx;\,f\in L^{1}(\R^{n})\cap L^{2}(\R^{n}).$$
%and
%$\Big(\|x\|=\Big(\displaystyle\sum_{i=1}^{n}x_{i}^{2}\Big)^{\frac{1}{2}}$
%(euclidean norm)\Big).
The Heisenberg-Pauli-Weyl inequality
for the Fourier-Bessel transform (one dimensional case) (see \cite{bowie,rosler}) is given by: for every $f \in L_\nu^2(\R_+)=L^2(\R_+,x^{2\nu+1}dx/2^\nu(\Gamma(\nu+1)))$ ($\nu>\frac{-1}{2}$) 
$$\Big\|\|x\|\,f\Big\|_{L_\nu^{2}(\R_+)}\Big\|\|y\|\,{\cal{F}}(f)\Big\|_{L_\nu^{2}(\R_+)}\geq\,
(\nu+1)\,\|f\|_{L_\nu^{2}(\R_+)}^{2},$$
where ${\cal{F}}(f)$ is the Fourier-Bessel defined by
$${\mathcal{F}}(f)(\xi)=\int_0^\infty f(x) j_{\nu}(\xi x)\frac{x^{2\nu+1}}{2^\nu\Gamma(\nu+1)}dx; \, \xi>0.$$
Here, $j_\nu$ is the
normalized Bessel function of the first kind and order $\nu>\frac{-1}{2}$ which is given by 
\begin{equation}\label{bfunction}
j_{\nu}(x)=\Gamma(\nu+1)\sum_{k=0}^{\infty}\frac{(-1)^{k}}{k!\,\Gamma(\nu+k+1)}(\frac{x}{2})^{2k}.
\end{equation}
Recently, the analogue of the Heisenberg-Pauli-Weyl inequality is established  for the Poly-axially operator (see \cite{moh ali}) and it is given by: 
for all $f\in L_{\alpha}^{2}(\R_{+}^{n})$
\begin{equation}\label{hpw}
\|f\|_{\alpha,2}^{2}\leq
\frac{2}{2|\alpha|+n}\,\Big\|\|x\|\,f\Big\|_{\alpha,2}\Big\|\|y\|\,{\mathcal{F}}_{\alpha}(f)\Big\|_{\alpha,2}.
\end{equation}
%where $A_\alpha=\dfrac{2}{2|\alpha|+n}$.
   \\
 \quad\quad In the present paper, we consider the Poly-axially operator (n-dimensional
Bessel operator) $\Delta_{\alpha}$ (see \cite{kc,CR}) defined for
$\alpha=(\alpha_{1},...,\alpha_{n})\in\R^{n}$ such that
$\alpha_{i}>-\dfrac{1}{2}$ for $i=1,...,n$ by
\begin{equation}\label{e1}
     \Delta_{\alpha}
     =\sum_{i=1}^{n}\frac{\p^{2}}{\p x_{i}^{2}}+\frac{2\alpha_{i}+1}{x_{i}}\frac{\p}{\p
     x_{i}}.
   \end{equation}
% It is also known as the Poly-axially operator (see \cite{kc}), and it is an important operator in analysis
%because of its applications in
%pure and applied Mathematics.\\
One can remark that:\\ $\bullet$ If $\alpha_{i}=-\dfrac{1}{2}$ for
$i=1,...,n$, then $\Delta_{\alpha}$ coincides with the Laplace
operator $\Delta_{n}$ on $\R^{n}$
$$\Delta_{\alpha}=\Delta_{n}=\sum_{i=1}^{n}\frac{\p^{2}}{\p x_{i}^{2}}.$$
$\bullet$ If $\alpha_{i}=-\dfrac{1}{2}$ for $i=1,...,n-1$ and
$\alpha_{n}>-\dfrac{1}{2}$, then $\Delta_{\alpha}$ coincides with
the Weinstein operator defined on $\R^{n}$(see
\cite{nb}, \cite{AS})
$$\Delta_{W}=\Delta_{n-1}+ L_{\alpha_{n}},$$
where $L_{\alpha_{n}}$ is the Bessel operator for the last variable
defined by
$$L_{\alpha_{n}}=\frac{\p^{2}}{\p x_{n}^{2}}+\frac{2\alpha_{n}+1}{x_{n}}\frac{\p}{\p x_{n}}.$$

%$(x)=\prod_{i=1}^{n}j_{\alpha_{i}}(x_{i});\quad x_{i}\in\R_{+}.$
%($j_{\gamma}$ is the normalized Bessel function
%of first kind and order $\gamma$.)\\
% More details are given in section 2.\\

Now, let m be a function in $L_{\alpha}^{2}(\R_{+}^{n})$ and let
$\sigma$ be a positive real number. The Poly-axially $L^{2}$-Multiplier
operators $T_{\alpha,m,\sigma}$ is defined (see \cite{CR}) for smooth functions
$\varphi$ on $\R_{+}^{n}$, by
%\begin{equation}\label{e18}
$$T_{\alpha,m,\sigma}\varphi(x)={\mathcal{F}}_{\alpha}^{-1}\Big(m_{\sigma}{\mathcal{F}}_{\alpha}(\varphi)\Big)(x),\;x\in\R_{+}^{n},$$
%\end{equation}
where the function $m_{\sigma}$ is given by $$m_{\sigma}(x)=m(\sigma
x),\;x\in\R_{+}^{n},$$

and ${\mathcal{F}}_{\alpha}$ is the n-dimensional Bessel Fourier transform
 defined for suitable function $f$ by
$${\mathcal{F}}_{\alpha}(f)(\lambda)=\int_{\R_{+}^{n}}f(x)\prod_{i=1}^{n} j_{\alpha_{i}}(\lambda_{i}x_{i})d\mu_{\alpha}
(x),$$ where $d\mu_{\alpha}(x)=c_{\alpha}\,\displaystyle{\prod_{i=1}^{n}}x_{i}^{2\alpha_{i}+1}dx_{i},$
 $c_{\alpha}
=\dfrac{1}{2^{|\alpha|}\prod_{i=1}^{n}\Gamma(\alpha_{i}+1)}$, ($|\alpha|=\alpha_{1}+...+\alpha_{n}$) and $j_{\alpha_{i}}$ is the
normalized Bessel function of the first kind and order $\alpha_i>\frac{-1}{2}$ given by $(\ref{bfunction})$.\\

 As well
%\begin{equation}\label{e19}
$$T_{\alpha,m,\sigma}\varphi(x)={\mathcal{F}}_{\alpha}^{-1}(m_{\sigma})*_{\alpha}\varphi(x),\;x\in\R_{+}^{n},$$
%\end{equation}
where
$${\mathcal{F}}_{\alpha}^{-1}(m_{\sigma})(x)=\frac{1}{\sigma^{2|\alpha|+n}}{\mathcal{F}}_{\alpha}^{-1}(m)(\frac{x}{\sigma}).$$
 The Poly-axially $L^{2}$-Multiplier
operators $T_{\alpha,m,\sigma}$ are a generalization of the usual
 linear multiplier operator $T_{m}$
% associated with a bounded function $m$ and given by $T_{m}(\varphi)={\mathcal{F}}^{-1}(m{\mathcal{F}}(\varphi))$,
 where $m$ is a bounded function and ${\mathcal{F}}(\varphi)$ denotes the ordinary transform on $\R^{n}$.\\
 % And many Mathematicians have been interested in these operators
 %as J.P.Anker \cite{JP}, J.J.Betancor \cite{JJ}, J.Gosselin, K.Stempak \cite{stem}, F.Soltani \cite{SF3}, A.Saoudi \cite{AS} and \cite{CR} .... \\
%The Heisenberg-Pauli-Weyl inequality for the Poly-axially operator is
%given by (A.M): for all $f\in L_{\alpha}^{2}(\R_{+}^{n})$

Our aim is to establish an analogue of the Heisenberg-Pauli-Weyl
uncertainty principle for these operators
. That is, we are going to show for $f\in
L_{\alpha}^{2}(\R_{+}^{n})$ that
$$\|f\|_{\alpha,2}\leq\frac{2}{2|\alpha|+n}\,
\Big\|\|y\|\,{\mathcal{F}}_{\alpha}(f)\Big\|_{\alpha,2}\Big\|\|x\|\,T_{\alpha,m,\sigma}f
\Big\|_{\sigma,\alpha,2};\, x,y\in\R_+^n,$$ provided $m \in
L_{\alpha}^{2}(\R_{+}^{n})\bigcap L_{\alpha}^{\infty}(\R_{+}^{n})$
satisfying the admissibility condition
$$\int_{0}^{\infty}|m_{\sigma}(x)|^{2}\frac{d\sigma}{\sigma}=1\,;\;x\in\R_{+}^{n}.$$
Then, based on the techniques of Donoho-Stark \cite{donoho}, we show
uncertainty principle of concentration type for the  Poly-axially $L^{2}$-Multiplier
operators $T_{\alpha,m,\sigma}$.\\ Let $E$ be a measurable subset of $\R^n$, $S$ be a measurable subset of $(0,\infty)\times\R^n$ and 
$f\in L_{\alpha}^{2}(\R_{+}^{n})$. If $f$ is
$\varepsilon$-concentrated on E and $T_{\alpha,m,\sigma}f$ is
$\delta$-concentrated on S, then
$$\|m\|_{\alpha,1}(\mu_{\alpha}(E))^{\frac{1}{2}}\Big(\int\int_{S}\frac{1}{\sigma^{2(2|\alpha|+n)}}d\Omega_{\alpha}(\sigma,x)\Big)^{\frac{1}{2}}
\geq1-(\varepsilon+\delta),$$
provided $m\in
L_{\alpha}^{1}(\R_{+}^{n})\bigcap L_{\alpha}^{2}(\R_{+}^{n})$
satisfying the precedent admissibility condition

\vskip0.2cm
 The contents of this work is the following. In section 2, we
recall some basic harmonic analysis results related to the Poly-axially
operator $\Delta_{\alpha}$ on $\R^{n}$. In section 3 and section 4,
we establish respectively the  Heisenberg-Pauli-Weyl uncertainty
principle and the Donoho-Stark’s uncertainty principle for the
Poly-axially $L_{\alpha}^{2}$-multiplier operators $T_{\alpha,m,\sigma}$.

\section{Harmonic analysis associated with the Poly-axially operator:}
\quad\quad In this section, we recall some basic results in harmonic
analysis related to the Poly-axially operator $\Delta_{\alpha}$
($\alpha=(\alpha_{1},...,\alpha_{n})\in\R^{n}$; $\alpha_{i}>\frac{-1}{2}$; $i=1,...,n$, and
$|\alpha|=\alpha_{1}+...+\alpha_{n}$).
Main references are \cite{kc,CR,KT,W}\\
In the following we denote by\\
$\bullet$ $\R_{+}^{n}=\{x=(x_{1},...,x_{n});x_{1}>0,...,x_{n}>0\}$.\\
A function $f$ on $\R^{n}$ is said to be even if it is even in each
variable.\\
 $\bullet$ ${\mathcal{C}}_{e}(\R^{n})$, the space of all
even
continuous functions on $\R^{n}$.\\
$\bullet$ ${\mathcal{C}}_{e,0}(\R^{n})$, the space of all even
continuous functions defined on $\R^{n}$ satisfying
$$\lim_{\|x\|\rightarrow\infty}f(x)=0\quad and\quad \|f\|_{{\mathcal{C}}_{e,0}}=\sup_{x\in\R_{+}^{n}}| f(x)|<+\infty,$$
where $\|x\|^{2}=x_{1}^{2}+...+x_{n}^{2}$.\\
$\bullet$ ${\mathcal{C}}_{e}^{m}(\R^{n})$, the space of all even $C^{m}$-functions on $\R^{n}$.\\
$\bullet$ ${\mathcal{S}}_{e}(\R^{n})$, the Schwartz space consists
of all even $C^{\infty}$-functions on $\R^{n}$ which are rapidly
decreasing as their derivatives, provided with the topology defined
for all $\gamma=(\gamma_{1},...,\gamma_{n})\in \N^{n}$, by
$$\rho_{m}(f)=\displaystyle \sup_{x\in\R^{n},\\ |\gamma|\leq m}(1+\|x\|^{2})^{m}|\p^{\gamma}f(x)|<\infty,\;\;\; \forall\;m\in\N, $$
with
$\p^{\gamma}f(x)=\dfrac{\p^{|\gamma|}f}{\p^{\gamma_{1}}x_{1}...\p^{\gamma_{n}}x_{n}}$.\\
$\bullet$ ${\mathcal{S}}'_{e}(\R^{n})$, the space of all even
tempered distributions in $\R^{n}$. That is the topological dual
 of ${\mathcal{S}}_{e}(\R^{n})$.\\
 $\bullet$
$L_{\alpha}^{p}(\R_{+}^{n},d\mu_{\alpha}(x))$, $1\leq p\leq\infty$,
the space of all measurable functions on $\R_{+}^{n}$ such that \quad
$\|f\|_{\alpha,p}=\Big(\displaystyle{\int}_{\R_{+}^{n}}| f(x)|^{p}d\mu_{\alpha}(x)\Big)^{\frac{1}{p}}<\infty,\quad 1\leq p<\infty,$
$$\|f\|_{\infty}=\displaystyle
ess\sup_{x\in\R_{+}^{n}}| f(x)|<\infty.$$ where\quad
$d\mu_{\alpha}(x)=c_{\alpha}\,\displaystyle\prod_{i=1}^{n}x_{i}^{2\alpha_{i}+1}dx_{i}$
\quad and\quad $c_{\alpha}
=\dfrac{1}{2^{|\alpha|}\prod_{i=1}^{n}\Gamma(\alpha_{i}+1)} .$\\
\\
 $\bullet$
$L_{\sigma,\alpha}^{p}(\R_{+}^{n},d\Omega_{\alpha}(\sigma,x))$,
$\sigma>0$, $1\leq p\leq\infty$, the space of all measurable
functions on $(0,\infty)\times\R_{+}^{n}$ such that \quad
$$\|f\|_{\sigma,\alpha,p}=\Big(\int_{\R_{+}^{n}}\int_{0}^{\infty}| f(\sigma,x)|^{p}d\Omega_{\alpha}(\sigma,x)\Big)^{\frac{1}{p}}<\infty,$$
where
$d\Omega_{\alpha}(\sigma,x):=\dfrac{d\sigma}{\sigma}d\mu_{\alpha}(x).$\\
\vskip0.5cm
 The function $j_{\alpha}$ is defined by $$
j_{\alpha}(x_{1},...,x_{n})=\prod_{i=1}^{n}j_{\alpha_{i}}(x_{i});\quad
(x_{1},...,x_{n})\in\R_{+}^{n},$$
where $j_{\alpha_{i}}$ is given by $(\ref{bfunction})$.
 %where $j_{\alpha_{i}}$ is the
%normalized Bessel function of the first kind and order $\alpha_i>\frac{-1}{2}$
%which is defined by
%$$j_{\alpha_{i}}(x_{i})=\Gamma(\alpha_{i}+1)\sum_{k=0}^{\infty}\frac{(-1)^{k}}{k!\,\Gamma(\alpha_{i}+k+1)}(\frac{x_{i}}{2})^{2k}.$$

$j_{\alpha}$ $(\alpha=(\alpha_{1},...,\alpha_{n}))$ verifies the following properties
\begin{equation}\label{j}
    | j_{\alpha}(x_{1},...,x_{n})|\leq 1\, ;\;
    (x_{1},...,x_{n})\;\in\R_{+}^{n},
\end{equation}
and  for $x=(x_{1},...,x_{n})$,  $ \in\R_{+}^{n}$, we have
$$\Delta_{\alpha}(j_{\alpha}(\lambda_{1} \,x_{1},...,\lambda_{n} \,x_{n}))=-\|\lambda\|^{2}j_{\alpha}(\lambda_{1} \,x_{1},...,\lambda_{n} \,x_{n});\;\;
\lambda=(\lambda_{1},...,\lambda_{n})\,\in\R_{+}^{n}.$$

\subsection{n-Dimensional Fourier-Bessel transform}
\begin{Def}
The Bessel Fourier transform is defined for suitable function f by
$${\mathcal{F}}_{\alpha}(f)(\lambda)=\int_{\R_{+}^{n}}f(x_{1},...,x_{n})j_{\alpha}(\lambda_{1} x_{1},...,\lambda_{n}x_{n})d\mu_{\alpha}(x_{1},...,x_{n});
\quad \lambda\in\R_{+}^{n}.$$
\end{Def}
It is shown in \cite{kc} and \cite{CR} that the Bessel
Fourier transform ${\mathcal{F}}_{\alpha}$ satisfies the following
properties.
\begin{Prop}\label{prop1}\quad
\begin{enumerate}
  \item For all $f\in L_{\alpha}^{1}(\R_{+}^{n})$,
  % the function${\mathcal{F}}_{\alpha}(f)$ is continous on $\R_{+}^{n}$ and
  we have
 $$ {\mathcal{F}}_{\alpha}(f)\in {\mathcal{C}}_{e,0}(\R_{+}^{n})\quad and \quad \|{\mathcal{F}}_{\alpha}(f)\|_{\alpha,\infty}\leq\|f\|_{\alpha,1}.$$
\item For $x\in\R_{+}^{n}$,\quad\quad
${\mathcal{F}}_{\alpha}(\Delta_{\alpha}(f))(x)=-\|x\|^{2}{\mathcal{F}}(f)(x).$
  \item The Bessel Fourier transform is a topological isomorphism from
  ${\mathcal{S}}_{e}(\R_{+}^{n})$ into itself. The inverse transform
  is given by
  \begin{equation}\label{e3}
   {\mathcal{F}}_{\alpha}^{-1}(f)(\lambda)={\mathcal{F}_{\alpha}}(f)(\lambda),\quad \lambda\in\R_{+}^{n},
   \end{equation}
  \end{enumerate}
\end{Prop}
\begin{Th}
(\underline{Inversion Formula})\\ Let f $\in
  L_{\alpha}^{1}(\R_{+}^{n})$ such that ${\mathcal{F}_{\alpha}}(f)\in
  L_{\alpha}^{1}(\R_{+}^{n})$, then we have
  \begin{equation}\label{e7}
  f(x)=\int_{\R_{+}^{n}}{\mathcal{F}}_{\alpha}(f)(\lambda)\prod_{i=1}^{n}j_{\alpha_{i}}(\lambda_{i}x_{i})d\mu_{\alpha}
  (\lambda);\quad x\in\R_{+}^{n}.
  \end{equation}
\end{Th}
\vskip0.5cm
The Bessel Fourier transform ${\mathcal{F}}_{\alpha}$ can be
extended to $L_{\alpha}^{2}(\R_{+}^{n})$ and we have
\begin{Th}
  (\underline{Plancherel Theorem})\\
The Fourier-Bessel transform ${\mathcal{F}_{\alpha}}$ is an isomorphism of $L_{\alpha}^{2}(\R_+^n)$ and for all $f\in L_{\alpha}^{2}(\R_+^n)$ we have
\begin{equation}\label{e6}
  \|{\mathcal{F}_{\alpha}}(f)\|_{\alpha,2}=\|f\|_{\alpha,2}.
  \end{equation}

\end{Th}

\subsection{The generalized translation operator associated with the Bessel operator}
\begin{Def}
%(See \cite{CHA},\cite{yil})
The generalized translation operator ${\cal{T}}_{x}^{\alpha}$,
$x\,\in\R_{+}^{n}$, associated with the Poly-axially operator
$\Delta_{\alpha}$ is defined on ${\mathcal{C}}_{e}(\R_{+}^{n})$, for
all y $\in\R_{+}^{n}$, by
$${\cal{T}}_{x}^{\alpha}f(y)=c'_{\alpha}\int_{[0,\pi]^{n}}f(X_{1},...,X_{n})\sin(\theta_{1})^{2\alpha_{1}}...\,\sin(\theta_{n})^{2\alpha_{n}}d\theta_{1}...d\theta_{n},$$
where \quad
$c'_{\alpha}=\prod_{i=1}^{n}\dfrac{\Gamma(\alpha_{i}+1)}{\sqrt{\pi}
\Gamma(\alpha_{i}+\frac{1}{2})}$\quad and\quad
$X_{i}=\sqrt{x_{i}^{2}+y_{i}^{2}-2x_{i}y_{i}\cos(\theta_{i})}$ ;
$i=1,...,n$.\\ As well, the generalized translation operator
${\cal{T}}_{x}^{\alpha}$ can be written as follow (see \cite{kc})
  \begin{equation}\label{e12}
  {\cal{T}}_{x}^{\alpha}f(y)=\int_{\R_{+}^{n}}\omega_{\alpha}(x,y,z)f(z)d\mu_{\alpha}(z),
  \end{equation}
  where\\
  $\omega_{\alpha}(x,y,z)=\\ \left\{
                            \begin{array}{ll}
                              \dfrac{c'_{\alpha}}{2^{2|\alpha|-n}}\dfrac{\prod_{i=1}^{n}\Big([z_{i}^{2}-(x_{i}^{2}-y_{i}^{2})][(x_{i}+y_{i})^{2}-z_{i}^{2}]\Big)
                              ^{\alpha_{i}-\frac{1}{2}}}{\prod_{i=1}^{n}(x_{i}y_{i}z_{i})^{2\alpha_{i}}},
                                \hbox{if z $\in\prod_{i=1}^{n}[| x_{i}-y_{i}|,x_{i}+y_{i}]$,} \\
                              0,  \hbox{otherwise.}
                            \end{array}
                           \right.$
\end{Def}
The following proposition summarizes some properties of the
translation operator:
\begin{Prop}\label{prop2}\quad
%The translation operator ${\cal{T}}_{x}^{\alpha}$, x $\in\R_{+}^{n}$,
%satisfies the following properties:
\begin{enumerate}
  \item For suitable function f, we have for all x,
  y $\in\R_{+}^{n}$
  %\begin{equation}\label{e8}
 $$ {\cal{T}}_{x}^{\alpha}f(y)={\cal{T}}_{y}^{\alpha}f(x)\quad and \quad {\cal{T}}_{0}^{\alpha}f(x)=f(x).$$
 % \end{equation}
  \item For all x, y $\in\R_{+}^{n}$ and $\lambda\in\C^{n}$, we have
  the following Product formula
 % \begin{equation}\label{e9}
 $$ {\cal{T}}_{x}^{\alpha}(j_{\alpha}(\lambda_{1}.,...,\lambda_{n}.))(y)=j_{\alpha}(\lambda_{1}x_{1},...,\lambda_{n}x_{n})j_{\alpha}(\lambda_{1}y_{1},...,
  \lambda_{n}y_{n}).$$
 % \end{equation}
  \item For f $\in L_{\alpha}^{p}(\R_{+}^{n})$, $1\leq p\leq\infty$,
  and x $\in\R_{+}^{n}$, then
  ${\cal{T}}_{x}^{\alpha}f\,\in L_{\alpha}^{p}(\R_{+}^{n})$
  and
  %\begin{equation}\label{e10}
 $$ \|{\cal{T}}_{x}^{\alpha}f\|_{\alpha,p}\leq\|f\|_{\alpha,p}.$$
  %\end{equation}

  \item For f,g $\in L_{\alpha}^{1}(\R_{+}^{n})$ and $x \in
  \R_{+}^{n}$, we have
 % \begin{equation}\label{et}
  $$\int_{\R_{+}^{n}}{\cal{T}}_{x}^{\alpha}f(y)g(y)d\mu_{\alpha}(y)=\int_{\R_{+}^{n}}f(y){\cal{T}}_{x}^{\alpha}g(y)d\mu_{\alpha}(y).$$
%\end{equation}
     \item For f $\in L_{\alpha}^{p}(\R_{+}^{n})$, $p=1,2$ and x
  $\in\R_{+}^{n}$, we have
  %\begin{equation}\label{e11}
  $${\mathcal{F}}_{\alpha}({\cal{T}}_{x}^{\alpha}f)(y)=\prod_{i=1}^{n}j_{\alpha_{i}}(y_{i}x_{i}){\mathcal{F}}_{\alpha}(f)(y); \quad y\,\in\R_{+}^{n}.$$
  %\end{equation}
 % \item For all f $\in {\mathcal{E}}_{*}(\R_{+}^{n})$ and y
 % $\in\R_{+}^{n}$, the function $x\longmapsto {\cal{T}}_{x}^{\alpha}f$
 % belongs to $ {\mathcal{E}}_{*}(\R_{+}^{n})$.
\end{enumerate}
\end{Prop}

\subsection{The generalized Bessel convolution product}
\begin{Def}
The generalized Bessel convolution product $*_{\alpha}$ for suitable
functions $f$ and $g$, is given by: for all $x\, \in\R_{+}^{n}$:
$$f*_{\alpha}g(x)=\int_{\R_{+}^{n}}{\cal{T}}_{x}^{\alpha}f(y)g(y)d\mu_{\alpha}(y).$$
\end{Def}
This convolution is commutative and associative, and it satisfies
the following properties:
\begin{Prop}\label{prop3}(Young's inequality)\quad
 Let p, q, r $\in[0,\infty]$, such that
  $\dfrac{1}{p}+\dfrac{1}{q}-\dfrac{1}{r}=1$ , then for all f $\in
  L_{\alpha}^{p}(\R_{+}^{n})$ and g $\in
  L_{\alpha}^{q}(\R_{+}^{n})$ the function $f*_{\alpha}g$ belongs to
  $ L_{\alpha}^{r}(\R_{+}^{n})$ and we have
  %\begin{equation}\label{e14}
  $$ \|f*_{\alpha}g\|_{\alpha,r}\leq \|f\|_{\alpha,p}\|g\|_{\alpha,q}.$$
  % \end{equation}
   \end{Prop}
   \begin{Prop}\quad
   \begin{enumerate}
  \item For all f, g $\in L_{\alpha}^{1} (\R_{+}^{n})$ (resp f, g $\in {\mathcal
  {S}}_{e}(\R_{+}^{n}$)) then
  $f*_{\alpha}g\,\in L_{\alpha}^{1}(\R_{+}^{n})$
  (resp $f*_{\alpha}g\,\in  {\mathcal
  {S}}_{e}(\R_{+}^{n})$) and we have
  %\begin{equation}\label{e13}
  $${\mathcal{F}}_{\alpha}(f*_{\alpha}g)= {\mathcal{F}}_{\alpha}(f){\mathcal{F}}_{\alpha}(g).$$
  %\end{equation}
  %\item Let $f\in L_{\alpha}^{1}(\R_{+}^{n})$ such that ${\mathcal{F}}_{\alpha}(f)\in L_{\alpha}^{1}(\R_{+}^{n})$, then for all $g\in
  %L_{\alpha}^{1}(\R_{+}^{n})$, we have $fg\in
 % L_{\alpha}^{1}(\R_{+}^{n})$ and $${\mathcal{F}}_{\alpha}(fg)={\mathcal{F}}_{\alpha}(f)*_{\alpha}{\mathcal{F}}_{\alpha}(g).$$
 \item For all  f, g $\in L_{\alpha}^{2}(\R_{+}^{n})$, we have
  %\begin{equation}\label{e17}
$${\mathcal{F}}_{\alpha}(fg)
    ={\mathcal{F}}_{\alpha}(f)*_{\alpha}{\mathcal{F}}_{\alpha}(g).$$
 
\end{enumerate}
\end{Prop}

\subsection{Poly-axially $L_{\alpha}^{2}$-multiplier operators on $\R_{+}^{n}$:}
These operators are studied in \cite{CR}, where the authors established some applications as Calder\'{o}n's reproducing formulas, best approximation formulas and extremal functions.
 \begin{Def}
Let m be a function in $L_{\alpha}^{2}(\R_{+}^{n})$ and let $\sigma$
be a positive real number. The Poly-axially  $L_{\alpha}^{2}$-multiplier operators
$T_{\alpha,m,\sigma}$ is defined for smooth functions $\varphi$ on
$\R_{+}^{n}$, by
\begin{equation}\label{e18}
T_{\alpha,m,\sigma}\varphi(x)={\mathcal{F}}_{\alpha}^{-1}\Big(m_{\sigma}{\mathcal{F}}_{\alpha}(\varphi)\Big)(x),\;x\in\R_{+}^{n},
\end{equation}
where $m_{\sigma}$ is the function given by
$$m_{\sigma}(x)=m(\sigma x).$$ As well
\begin{equation}\label{e19}
T_{\alpha,m,\sigma}\varphi(x)={\mathcal{F}}_{\alpha}^{-1}(m_{\sigma})*_{\alpha}\varphi(x),\;x\in\R_{+}^{n},
\end{equation}
where
$${\mathcal{F}}_{\alpha}^{-1}(m_{\sigma})(x)=\frac{1}{\sigma^{2|\alpha|+n}}{\mathcal{F}}_{\alpha}^{-1}(m)(\frac{x}{\sigma}).$$
% and $|\alpha|=\alpha_{1}+...+\alpha_{n}.$
 \end{Def}
 
\begin{Lem}\label{lem1}
Let $m,\varphi\in L_{\alpha}^{1}(\R_{+}^{n})\bigcap
L_{\alpha}^{2}(\R_{+}^{n})$, then
\begin{equation}\label{e22}
    T_{\alpha,m,\sigma}\varphi(\lambda)=\frac{1}{\sigma^{2|\alpha|+n}}\int_{\R_{+}^{n}}\Theta_{\alpha}(y,
    \lambda,m)\,\varphi(y)\,d\mu_{\alpha}(y),
\end{equation}
where
$\Theta_{\alpha}(y,\lambda,m)=\displaystyle \int_{\R_{+}^{n}}m(x)\prod_{i=1}^{n}j_{\alpha_{i}}(x_i\dfrac{\lambda_{i}}{\sigma})\prod_{i'=1}^{n}j_{\alpha_{i'}}(x_i\dfrac{y_{i'}}{\sigma})
d\mu_{\alpha}(x)$.
\end{Lem}
{\bf Proof:}\\
Using the relations ($\ref{e18}$), ($\ref{e3}$) and Fubini's
theorem, we have
\begin{eqnarray*}
% \nonumber to remove numbering (before each equation)
  T_{\alpha,m,\sigma}\varphi(\lambda)&=&\int_{\R_{+}^{n}}m(\sigma x){\mathcal{F}}_{\alpha}(\varphi)(x)
  \prod_{i=1}^{n}j_{\alpha_{i}}(\lambda_{i}x_{i})
  d\mu_{\alpha}(x)  \\
  % &=& \int_{\R_{+}^{n}}m(\sigma x)\Big(\int_{\R_{+}^{n}}\varphi(y)\prod_{i=1}^{n}j_{\alpha_{i}}(y_{i}x_{i})
  %d\mu_{\alpha}(y)\Big)\prod_{i=1}^{n}j_{\alpha_{i}}(\lambda_{i}x_{i})
  %d\mu_{\alpha}(x) \\
   &=&\int_{\R_{+}^{n}}\Big(\int_{\R_{+}^{n}}m(\sigma x)\prod_{i=1}^{n}j_{\alpha_{i}}(\lambda_{i}x_{i})\prod_{i'=1}^{n}j_{\alpha_{i'}}(y_{i'}x_{i})d\mu_{\alpha}(x)\Big)
  \varphi(y)d\mu_{\alpha}(y).
\end{eqnarray*}
Next, with a simple change of variable, we get
$$ T_{\alpha,m,\sigma}\varphi(\lambda)=\frac{1}{\sigma^{2|\alpha|+n}}\int_{\R_{+}^{n}}\Big(\int_{\R_{+}^{n}}m(x)
  \prod_{i=1}^{n}j_{\alpha_{i}}(\frac{\lambda_{i}}{\sigma}\,x_{i})\prod_{i'=1}^{n}j_{\alpha_{i'}}
(\frac{y_{i'}}{\sigma}\,x_{i})d\mu_{\alpha}(x)\Big)
  \varphi(y)d\mu_{\alpha}(y),$$
which completes the proof.\\
%More properties and details are given in \cite{CR}

%\begin{Th}\label{th1}
%Let m be a function in $L_{\alpha}^{2}(\R_{+}^{n})$ satisfying the
%admissibility condition
%\begin{equation}\label{ac}
%\int_{0}^{\infty}|m_{\sigma}(x)|^{2}\frac{d\sigma}{\sigma}=1\,;\;x\in\R_{+}^{n}
%\end{equation}
  %Then, for all $\varphi\in L_{\alpha}^{2}(\R_{+}^{n})$, we
  %have
  %\begin{equation}\label{e25}
%\int_{\R_{+}^{n}}\int_{0}^{\infty}
%\|T_{\alpha,m,\sigma}\varphi\|_{\sigma,\alpha,2}=\|\varphi\|_{H_{\alpha}^{0}}=\|\varphi\|_{\alpha,2}.
%  \end{equation}
%\end{Th}

\section{Heisenberg-Pauli-Weyl uncertainty principle}
\quad\quad This section is devoted to establish
Heisenberg-Pauli-Weyl uncertainty principle for the Poly-axially $L^{2}$-
multiplier operators $T_{\alpha,m,\sigma}$.
\begin{Th}\label{unt1}
Let m be a function in $L_{\alpha}^{2}(\R_{+}^{n})\bigcap
L_{\alpha}^{\infty}(\R_{+}^{n})$ satisfying the following
admissibility condition
\begin{equation}\label{ac}
\int_{0}^{\infty}|m_{\sigma}(x)|^{2}\frac{d\sigma}{\sigma}=1\,;\;x\in\R_{+}^{n}.
\end{equation}
 Then, for $f\in
L_{\alpha}^{2}(\R_{+}^{n})$, we have
\begin{equation}
\|f\|_{\alpha,2}\leq \frac{2}{2|\alpha|+n}\,
\Big\|\|y\|\,{\mathcal{F}}_{\alpha}(f)\Big\|_{\alpha,2}\Big\|\|x\|\,T_{\alpha,m,\sigma}f
\Big\|_{\sigma,\alpha,2}.
\end{equation}
 %(C depends on n???)
\end{Th}
{\bf Proof:}\\ Let $f\in L_{\alpha}^{2}(\R_{+}^{n})$. Assume that
$\Big\|\|y\|\,{\mathcal{F}}_{\alpha}(f)\Big\|_{\alpha,2}<\infty$\quad
and\quad $\Big\|\|x\|\,T_{\alpha,m,\sigma}f
\Big\|_{\sigma,\alpha,2}<\infty.$\\ According to inequality
($\ref{hpw})$ and relation ($\ref{e18}$), we get
\begin{eqnarray*}
% \nonumber to remove numbering (before each equation)
  \int_{\R_{+}^{n}}|T_{\alpha,m,\sigma}f(x)|^{2}d\mu_{\alpha}(x) %&\leq&\frac{2}{2|\alpha|+n}\,\Big(\int_{\R_{+}^{n}}\|x\|^{2}|T_{\alpha,m,\sigma}f(x)|^{2}
 % d\mu_{\alpha}(x)\Big)^{\frac{1}{2}}\\
 % &\times&\Big(\int_{\R_{+}^{n}}\|y\|^{2}|{\mathcal{F}}_{\alpha}(T_{\alpha,m,\sigma}f(.))(y)|^{2}
  %d\mu_{\alpha}(y)\Big)^{\frac{1}{2}}   \\
   &\leq&\frac{2}{2|\alpha|+n}\,\Big(\int_{\R_{+}^{n}}\|x\|^{2}|T_{\alpha,m,\sigma}f(x)|^{2}
  d\mu_{\alpha}(x)\Big)^{\frac{1}{2}}\\
  &\times&\Big(\int_{\R_{+}^{n}}\|y\|^{2}|m(\sigma y){\mathcal{F}}_{\alpha}(f)(y))|^{2}
  d\mu_{\alpha}(y)\Big)^{\frac{1}{2}}
\end{eqnarray*}
Integrating over $[0,\infty[$ with respect to the measure
$\dfrac{d\sigma}{\sigma}$, we obtain
\begin{eqnarray*}
\int_{0}^{\infty}\Big(\int_{\R_{+}^{n}}|T_{\alpha,m,\sigma}f(x)|^{2}d\mu_{\alpha}(x)\Big)\frac{d\sigma}{\sigma}&\leq&
\frac{2}{2|\alpha|+n}\,\int_{0}^{\infty}\Big(\int_{\R_{+}^{n}}\|x\|^{2}|T_{\alpha,m,\sigma}f(x)|^{2}
  d\mu_{\alpha}(x)\Big)^{\frac{1}{2}}\\
  &\times&\Big(\int_{\R_{+}^{n}}\|y\|^{2}|m(\sigma
y){\mathcal{F}}_{\alpha}(f)(y))|^{2}
  d\mu_{\alpha}(y)\Big)^{\frac{1}{2}}\frac{d\sigma}{\sigma}.
\end{eqnarray*}
Using (\cite{CR}Theorem 2.4 (i)) and Schwarz's inequality, we have
\begin{eqnarray*}
\|f\|_{\alpha,2}&\leq&\frac{2}{2|\alpha|+n}\,\Big(\int_{0}^{\infty}\int_{\R_{+}^{n}}\|x\|^{2}|T_{\alpha,m,\sigma}f(x)|^{2}
  d\mu_{\alpha}(x)\frac{d\sigma}{\sigma}\Big)^{\frac{1}{2}}\\
  &\times&\Big(\int_{0}^{\infty}\int_{\R_{+}^{n}}\|y\|^{2}|m(\sigma
y){\mathcal{F}}_{\alpha}(f)(y))|^{2}
  d\mu_{\alpha}(y)\frac{d\sigma}{\sigma}\Big)^{\frac{1}{2}}.
\end{eqnarray*}
Now, by Fubini-Tonneli's theorem and the admissibility condition
($\ref{ac}$), we get
\begin{eqnarray*}
\|f\|_{\alpha,2}&\leq&\frac{2}{2|\alpha|+n}\,\Big(\int_{0}^{\infty}\int_{\R_{+}^{n}}\|x\|^{2}|T_{\alpha,m,\sigma}f(x)|^{2}
  d\mu_{\alpha}(x)\frac{d\sigma}{\sigma}\Big)^{\frac{1}{2}}\\
  &\times&\Big(\int_{\R_{+}^{n}}\|y\|^{2}|{\mathcal{F}}_{\alpha}(f)(y))|^{2}
  d\mu_{\alpha}(y)\Big)^{\frac{1}{2}}.
\end{eqnarray*}
This completes the proof of the theorem.

\begin{Th}
Let m be a function in $L_{\alpha}^{2}(\R_{+}^{n})\bigcap
L_{\alpha}^{\infty}(\R_{+}^{n})$ satisfying the admissibility
condition ($\ref{ac}$) and $a,\,b\in[1,\infty)$. Let
$\epsilon\in\R$ such that
$a\epsilon=(1-\epsilon)b$, then for $f\in
L_{\alpha}^{2}(\R_{+}^{n})$, we have
\begin{equation}
\|f\|_{\alpha,2}\leq \Big(\frac{2}{2|\alpha|+n}\Big)^{a\epsilon}\,
\Big\|\|x\|^{a}\,T_{\alpha,m,\sigma}f\Big\|_{\sigma,\alpha,2}^{\epsilon}\,\Big\|\|y\|^{b}\,{\mathcal{F}}_{\alpha}(f)\Big\|_{\alpha,2}
^{1-\epsilon}.
\end{equation}
\end{Th}
{\bf Proof:}\\ Let $f\in L_{\alpha}^{2}(\R_{+}^{n})$, $f\neq 0$. Assume that
$\Big\|\|x\|^{a}\,T_{\alpha,m,\sigma}f\Big\|_{\sigma,\alpha,2}^{\epsilon}<\infty$
and $\Big\|\|y\|^{b}\,{\mathcal{F}}_{\alpha}(f)\Big\|_{\alpha,2}
^{1-\epsilon}<\infty.$ Now, for all $a\geq 1$, we have
$$\Big\|\|x\|^{a}\,T_{\alpha,m,\sigma}f\Big\|_{\sigma,\alpha,2}^{\frac{1}{a}}\,\Big\|T_{\alpha,m,
\sigma}f\Big\|_{\sigma,\alpha,2}^{\frac{1}{a'}}=\Big\|\|x\|^{2}\,T_{\alpha,m,\sigma}(f)^{\frac{2}{a}}\Big\|_{\sigma,\alpha,2}^{\frac{1}{2}}\,\Big\|
T_{\alpha,m,\sigma}(f)^{\frac{2}{a'}}\Big\|_{\sigma,\alpha,2}^{\frac{1}{2}},$$
where $\dfrac{1}{a}+\dfrac{1}{a'}=1.$ 
 Using H\"{o}lder's inequality, we get
$$\Big\|\|x\|\,T_{\alpha,m,\sigma}f\Big\|_{\sigma,\alpha,2}\leq\Big\|\|x\|^{a}\,T_{\alpha,m,\sigma}f\Big\|_{\sigma,\alpha,2}^{\frac{1}{a}}\,\Big\|T_{\alpha,m,
\sigma}f\Big\|_{\sigma,\alpha,2}^{\frac{1}{a'}} .$$
%\begin{eqnarray*}
% \nonumber to remove numbering (before each equation)
 % \Big\|\|x\|\,T_{\alpha,m,\sigma}f\Big\|_{\alpha,2} &=&\Big\|\|x\|^{2}\,T_{\alpha,m,\sigma}(f)^{2}\Big\|_{\alpha,1}^{\frac{1}{2}}  \\
 %  &\leq&\Big\|\|x\|^{2}\,T_{\alpha,m,\sigma}(f)^{\frac{2}{\beta}}\Big\|_{\alpha,\beta}^{\frac{1}{2}}\,\Big\|T_{\alpha,m,\sigma}(f)
 %  ^{\frac{2}{\beta'}}\Big\|_{\alpha,\beta'}^{\frac{1}{2}}   \\
 %  &\leq&\Big\|\|x\|^{\beta}\,T_{\alpha,m,\sigma}f\Big\|_{\alpha,2}^{\frac{1}{\beta}}\,\Big\|{\mathcal{T}}_{\alpha,m,%
%\sigma}f\Big\|_{\alpha,2}^{\frac{1}{\beta'}}\\
%   &\leq&\Big\|\|x\|^{\beta}\,T_{\alpha,m,\sigma}f\Big\|_{\alpha,2}^{\frac{1}{\beta}}\,\Big\|m\Big\|_{\alpha,\infty}\,\Big\|f\Big\|_{\alpha,2}.
%\end{eqnarray*}
Referring to (\cite{CR}, Theorem 2.4), we have for all $ a\geq1$
 \begin{equation}\label{e20}
\Big\|\|x\|\,T_{\alpha,m,\sigma}f\Big\|_{\sigma,\alpha,2}
\leq\Big\|\|x\|^{a}\,T_{\alpha,m,\sigma}f\Big\|_{\sigma,\alpha,2}^{\frac{1}{a}}\,\Big\|f\Big\|_{\alpha,2}^{\frac{1}{a'}},
 \end{equation}
 with equality if $a=1$.\\
Using the same method, for all $b\geq1$ and $\dfrac{1}{b}+\dfrac{1}{b'}=1$, we have
\begin{equation}\label{e21}
    \Big\|\|y\|\,{\mathcal{F}}_{\alpha}(f)\Big\|_{\alpha,2}\leq\Big\|\|y\|^{b}\,{\mathcal{F}}_{\alpha}(f)\Big\|_{\alpha,2}^{\frac{1}{b}}\,\Big\|f\Big\|_
    {\alpha,2}^{\frac{1}{b'}},
\end{equation}
with equality if $b=1$.\\
 According to relations ($\ref{e20}$) and ($\ref{e21}$) and by the fact that
 $a\epsilon=(1-\epsilon)b$, we get
 $$\Big(\frac{\Big\|\|x\|\,T_{\alpha,m,\sigma}f\Big\|_{\sigma,\alpha,2}\Big\|\|y\|\,{\mathcal{F}}_{\alpha}(f)\Big\|_{\alpha,2}}
 {\Big\|f\Big\|_{\alpha,2}^{\frac{1}{a'}}\Big\|f\Big\|_
    {\alpha,2}^{\frac{1}{b'}}}\Big)^{a\epsilon} \leq  \Big\|\|x\|^{a}\,T_{\alpha,m,\sigma}f\Big\|_{\sigma,\alpha,2}^{\epsilon}
   \Big\|\|y\|^{b}\,{\mathcal{F}}_{\alpha}(f)\Big\|_{\alpha,2}^{1-\epsilon},$$
with equality if $a=b=1$.\\
Now, using Theorem $\ref{unt1}$, we have
$$\|f\|_{\alpha,2}\leq \Big(\frac{2}{2|\alpha|+n}\Big)^{a\epsilon}\,
\Big\|\|x\|^{a}\,T_{\alpha,m,\sigma}f\Big\|_{\sigma,\alpha,2}^{\epsilon}\,\Big\|\|y\|^{b}\,{\mathcal{F}}_{\alpha}(f)\Big\|_{\alpha,2}
^{1-\epsilon}.$$

\section{Donoho-Stark's uncertainty principle}
\begin{Def} (see \cite{fs,KH})\quad
\begin{enumerate}
\item[(i)] Let E be a measurable subset of $\R_{+}^{n}$, we say that the function $f\in L_{\alpha}^{2}(\R_{+}^{n})$ is
$\varepsilon$-concentrated on E if
\begin{equation}\label{e22}
\|f-\chi_{E}f\|_{\alpha,2}\leq \varepsilon\|f\|_{\alpha,2},
\end{equation}
where $\chi_{E}$ is the indicator function of the set E.
\item[(ii)]  Let S be a measurable subset of $(0,\infty)\times\R_{+}^{n}$ and let $f\in
L_{\alpha}^{2}(\R_{+}^{n})$. We say that $T_{\alpha,m,\sigma}f$ is a
$\delta$-concentrated on S if
\begin{equation}\label{e23}
\|T_{\alpha,m,\sigma}f-\chi_{S}T_{\alpha,m,\sigma}f\|_{\sigma,\alpha,2}\leq
\delta\|T_{\alpha,m,\sigma} f\|_{\sigma,\alpha,2}.
\end{equation}
 Where $\chi_{S}$ is the indicator function of the set S.\\
% $$\|f\|_{\sigma,\alpha,2}=\Big(\int_{\R_{+}^{n}}\int_{+}^{\infty}|f(\sigma,x)|^{2}d\Omega_{\alpha}(\sigma,x)\Big)^{\frac{1}{2}}$$
% and  $\Omega_{\alpha}$ is the measure on $(0,\infty)\times
%\R_{+}^{n}$ given by
%$$d\Omega_{\alpha}(\sigma,x):=\frac{d\sigma}{\sigma}d\mu_{\alpha}(x).$$
\end{enumerate}
\end{Def}

\begin{Th}
Let $f$ be a function in $L_{\alpha}^{2}(\R_{+}^{n})$ and $m\in
L_{\alpha}^{1}\bigcap L_{\alpha}^{2}(\R_{+}^{n})$ satisfying the
admissibility condition ($\ref{ac}$). If $f$ is
$\varepsilon$-concentrated on E and $T_{\alpha,m,\sigma}f$ is
$\delta$-concentrated on S, then
\begin{equation}
\|m\|_{\alpha,1}(\mu_{\alpha}(E))^{\frac{1}{2}}\Big(\int\int_{S}\frac{1}{\sigma^{2(2|\alpha|+n)}}d\Omega_{\alpha}(\sigma,x)\Big)^{\frac{1}{2}}
\geq1-(\varepsilon+\delta).
\end{equation}
%where $\Omega_{\alpha}$ is the measure on $(0,\infty)\times
%\R_{+}^{n}$ given by
%$$d\Omega_{\alpha}(\sigma,x):=\frac{d\sigma}{\sigma}d\mu_{\alpha}(x).$$
\end{Th}

{\bf Proof:}\\
Let $f\in L_{\alpha}^{2}(\R_{+}^{n})$ and $m\in
L_{\alpha}^{1}(\R_{+}^{n})\bigcap L_{\alpha}^{2}(\R_{+}^{n})$
satisfying ($\ref{ac}$). Assume that $\mu_{\alpha}(E)<\infty$ and
$\dint\dint_{S}\dfrac{1}{\sigma^{2(2|\alpha|+n)}}d\Omega_{\alpha}(\sigma,x)<\infty.$\\
According to relations ($\ref{e22}$) and ($\ref{e23}$), we have
\begin{eqnarray*}
% \nonumber to remove numbering (before each equation)
   \|T_{\alpha,m,\sigma}f-\chi_{S}T_{\alpha,m,\sigma}(\chi_{E}f)\|_{\sigma,\alpha,2}%&=&\|T_{\alpha,m,\sigma}
  % f-(\chi_{S}T_{\alpha,m,\sigma}(f-(f-\chi_{E}f))\|_{\sigma,\alpha,2}  \\
  % &=& \|T_{\alpha,m,\sigma}
  % f-(\chi_{S}T_{\alpha,m,\sigma}f-\chi_{S}T_{\alpha,m,\sigma}(f-\chi_{E}f)\|_{\sigma,\alpha,2} \\
   &\leq& \|{\mathcal{T}}_{\sigma,\alpha,m}
   f-\chi_{S}T_{\alpha,m,\sigma}f\|_{\sigma,\alpha,2}+\|\chi_{S}T_{\alpha,m,\sigma}(f-\chi_{E}f)\|_{\sigma,\alpha,2} \\
   &\leq&\delta\|T_{\alpha,m,\sigma}f\|_{\sigma,\alpha,2}+\|T_{\alpha,m,\sigma}(f-\chi_{E}f)\|_{\sigma,\alpha,2}.
\end{eqnarray*}
And by (\cite{CR}Theorem 2.4 (i)), we obtain
\begin{eqnarray*}
  \|T_{\alpha,m,\sigma}f-\chi_{S}T_{\alpha,m,\sigma}(\chi_{E}f)\|_{\sigma,\alpha,2} &\leq&\delta\|f\|_{\alpha,2}+\|f-\chi_{E}f\|_{\alpha,2}  \\
   &\leq&(\delta+\varepsilon)\|f\|_{\alpha,2}.
\end{eqnarray*}
Then, using the triangle inequality we have
\begin{eqnarray*}
% \nonumber to remove numbering (before each equation)
  \|T_{\alpha,m,\sigma}f\|_{\sigma,\alpha,2} &\leq& \|\chi_{S}T_{\alpha,m,\sigma}(\chi_{E}f)\|_{\sigma,\alpha,2}+
  \|T_{\alpha,m,\sigma}f-\chi_{S}T_{\alpha,m,\sigma}(\chi_{E}f)\|_{\sigma,\alpha,2} \\
   &\leq&\|\chi_{S}T_{\alpha,m,\sigma}(\chi_{E}f)\|_{\sigma,\alpha,2}+(\delta+\varepsilon)\|f\|_{\alpha,2}.
\end{eqnarray*}
Next, since $m$, $\chi_{E}f\in L_{\alpha}^{1}(\R_{+}^{n})\bigcap
L_{\alpha}^{2}(\R_{+}^{n})$, thus,  using Lemma $\ref{lem1}$, we have

\begin{eqnarray*}
% \nonumber to remove numbering (before each equation)
  |T_{\alpha,m,\sigma}(\chi_{E}f)(x)| %&=& |{\mathcal{F}}_{\alpha}^{-1}(m_{\sigma}{\mathcal{F}}_{\alpha}(\chi_{E}f))| \\
  % &\leq& \|{\mathcal{F}}_{\alpha}^{-1}(m_{\sigma}{\mathcal{F}}_{\alpha}(\chi_{E}f))\|_{\alpha,\infty}  \\
  % &\leq& \|m_{\sigma}{\mathcal{F}}_{\alpha}(\chi_{E}f)\|_{\alpha,1} \\
  % &\leq&\|m_{\sigma}\|_{\alpha,1}\,\|{\mathcal{F}}_{\alpha}(\chi_{E}f)\|_{\alpha,\infty}  \\
   &\leq& \frac{1}{\sigma^{2|\alpha|+n}}\|m\|_{\alpha,1}\,\|\chi_{E}f\|_{\alpha,1}\\
   &\leq&
   \frac{1}{\sigma^{2|\alpha|+n}}\|m\|_{\alpha,1}\,\|f\|_{\alpha,2}\,(\mu_{\alpha}(E))^{\frac{1}{2}}.
\end{eqnarray*}
Then, by the fact that
$\|\chi_{S}T_{\alpha,m,\sigma}(\chi_{E}f)\|_{\sigma,\alpha,2}=\Big(\dint\dint_{S}|T_{\alpha,m,\sigma}(\chi_{E}f)(x)|
^{2}d\Omega_{\alpha}(\sigma,x)\Big)^{\frac{1}{2}}$, we get
$$\|\chi_{S}T_{\alpha,m,\sigma}(\chi_{E}f)\|_{\sigma,\alpha,2}\leq\|m\|_{\alpha,1}\,\|f\|_{\alpha,2}\,(\mu_{\alpha}(E))^{\frac{1}{2}}
\Big(\int\int_{S}\dfrac{d\Omega_{\alpha}(\sigma,x)}{\sigma^{2(2|\alpha|+n)}}\Big)^{\frac{1}{2}}.$$
Therefore
$$\|T_{\alpha,m,\sigma}f\|_{\sigma,\alpha,2}\leq\|m\|_{\alpha,1}\|f\|_{\alpha,2}(\mu_{\alpha}(E))^{\frac{1}{2}}
\Big(\int\int_{S}\dfrac{d\Omega_{\alpha}(\sigma,x)}{\sigma^{2(2|\alpha|+n)}}\Big)^{\frac{1}{2}}+(\delta+\varepsilon)\|f\|_{\alpha,2}.$$
By (\cite{CR},theorem 2.4 (i)), it follows
$$\|f\|_{\alpha,2}\leq\|m\|_{\alpha,1}\|f\|_{\alpha,2}(\mu_{\alpha}(E))^{\frac{1}{2}}
\Big(\int\int_{S}\dfrac{d\Omega_{\alpha}(\sigma,x)}{\sigma^{2(2|\alpha|+n)}}\Big)^{\frac{1}{2}}+(\delta+\varepsilon)\|f\|_{\alpha,2}.$$
Then
$$\|m\|_{\alpha,1}(\mu_{\alpha}(E))^{\frac{1}{2}}\Big(\int\int_{S}\frac{1}{\sigma^{2(2|\alpha|+n)}}d\Omega_{\alpha}(\sigma,x)\Big)^{\frac{1}{2}}
\geq1-(\varepsilon+\delta),$$ which is the desired result.

\begin{Rem}
If
$S\subset\{(\sigma,x)\in(0,\infty)\times\R_{+}^{n}:\,\sigma\geq\xi\}$
for some $\xi>0$, one suppose that
$$\varrho=\max\{\frac{1}{\sigma}:\,(\sigma,x)\in S\, for\, some \,x\in\R_{+}^{n}\}. $$
Then by the previous theorem, we conclude that
$$\varrho^{2|\alpha|+n}\|m\|_{\alpha,1}(\mu_{\alpha}(E))^{\frac{1}{2}}(\Omega_{\alpha}(S))^{\frac{1}{2}}\geq 1-(\varepsilon+\delta).$$

\end{Rem}

\end{document}